\magnification=1200
\font\big=cmbx10 scaled\magstep2
\font\sc=cmcsc10
\font\ib=cmmib10
\font\tenmsym=msbm10
\font\eightmsym=msbm8
\font\teneufm=eufm10
\textfont 9=\tenmsym \scriptfont 9=\eightmsym 
\def\bb {\fam9 } 
\textfont 8=\teneufm \scriptfont 8=\teneufm

\centerline{\big Matrix Vieta Theorem}\bigskip
 
\centerline{{\sc Dmitry FUCHS and Albert SCHWARZ}\footnote*{This work was 
partially supported by NSF grant DMS-9201366. Research at MSRI was supported 
by NSF grant DMS-9022140.}}
 
\centerline{Department of Mathematics}
 
\centerline{University of California}
 
\centerline{Davis Ca 95616, USA}\bigskip\bigskip
 
\centerline{\bf 1. Introduction}\medskip
 
Consider a matrix algebraic equation$$X^n+A_1X^{n-1}+\dots +A_n=0,\eqno(1)$$ 
where the coefficients $A_1,\dots ,A_n$ as well as solutions $X$ are supposed 
to be square complex matrices of some order $k$. For an usual algebraic 
equation of degree $n$ the classical Vieta formulas express the coefficients 
in terms of the $n$ solutions. However, a matrix degree $n$ equation 
generically has $\displaystyle{nk\choose k}$ rather  than $n$ solutions. 
(Throughout this article the word {\it generic} refers to a Zariski open set.) 
 
We call $n$ solutions $X_1,\dots ,X_n$ of the equation (1) {\it independent} 
if they determine the coefficients $A_1,\dots ,A_n$ (a more technical 
explanation of independence see in Section 2). However, the expressions of 
$A_1,\dots ,A_n$ in terms of $X_1,\dots ,X_n$ are much less elegant, than 
Vieta formulas (see, e.g. formulas (3) -- (5) below). Still there are some 
relations between $A$'s and $X$'s which are very similar to Vieta.\smallskip 
 
{\sc Theorem} 1.1. {\it If solutions $X_1,\dots ,X_n$ of the equation $(1)$ 
are independent, then$$\eqalign{{\rm tr}\, A_1&=-({\rm tr}\, X_1+\dots +{\rm 
tr}\, X_n),\cr {\rm det}\, A_n&=(-1)^{nk}{\rm det}\, X_1\ldots{\rm det}\, 
X_n.\cr}\eqno(2)$$}\smallskip 
 
Theorem 1.1 is proved in Section 4 (a more direct proof for $n=2$ is given in 
Section 3). In Section 5 we discuss a generalization of Theorem 1.1 from 
complex matrix algebras to arbitrary associative unitary rings.\smallskip 
 
Theorem 1.1 is elementary, however it is connected with some constructions of 
modern Mathematics. For every associative algebra $A$ one can construct a 
linear space $F(A)=A^2/[A,A]$. This space appeared in [K] as the space of 0-
forms on a non-commutative formal manifold. (Such a manifold is determined by 
a free associative algebra $A$. For the free associative algebra $A$ with 
generators $a_1,\ldots ,a_n$ the space $F(A)$ is spanned by cyclic words of 
length $\ge2$ with letters $a_1,\ldots a_n$.) The space $F(A)$ appeared also 
in [GS] and [AW] in relation to other problems. It is easy to understand that 
the formula for ${\rm tr}\, A_1$ can be interpreted as a non-trivial identity 
in $F(A)$.\hfill\eject 
 
\centerline{\bf 2. Independent matrices.}\medskip
 
{\sc Definition} 2.1. Matrices $X_1,\dots ,X_n$ are called {\it independent} 
if the bloc Vandermonde determinant is not zero:$$\left|\matrix{I&I&\dots&I\cr 
X_1&X_2&\dots&X_n\cr \multispan4\quad$\cdots\cdots\cdots\cdots 
\cdots\cdots\cdots$\cr X_1^{n-1}&X_2^{n-1}&\dots&X_n^{n-1}}\right|\ne0.$$ 
 
For $n=2$ the independence condition means that ${\rm det}\, (X_1-X_2)\ne0$. 
For $n\ge3$ it does not imply and is not implied by the condition ${\rm 
det}\,(X_i-X_j)\ne0,\, 1\le i<j\le n$.
 
It is obvious that a generic equation (1) has $n$ independent solutions 
(otherwise the above determinant would have vanished for any $n$ solutions of 
any equation (1)). It is clear also that the matrices $X_1,\dots ,X_n$ are 
independent if and only is there exist unique $A_1,\dots ,A_n$ such that 
$X_1,\dots ,X_n$ satisfy the equation (1). In other words, for independent 
$X_1,\dots ,X_n$ the matrices $A_1,\dots ,A_n$ may be expressed via $X_1,\dots 
,X_n$. For example, if $n=2$ 
then$$\eqalign{A_1&=-(X_1^2-X_2^2)(X_1-X_2)^{-1}, 
\cr A_2&=-X_1^2+(X_1^2-X_2^2)(X_1-X_2)^{-1}X_1.\cr}\eqno(3)$$For $n\ge3$ it is 
impossible to write a formula valid for all independent matrices; for example, 
if $n=3$, then$$\eqalign{A_1=-((X_1^3-X_2^3)(X_1&-X_2)^{-1}-(X_1^3-X_3^3)(X_1-
X_3)^{-1})\cr &((X_1^2-X_2^2)(X_1-X_2)^{-1}-(X_1^2-X_3^2)(X_1-X_3)^{-1})^{-
1}\cr}\eqno(4)$$ provided that the right hand side exists (which does not 
follow from $X_1,X_2,X_3$ being independent); otherwise the expression will be 
different.
 
An expression of $A_1,\ldots ,A_n$ as functions of $X_1,\ldots ,X_n$, valid 
for generic independent $X_1,\ldots ,X_n$, may be given in terms of 
Gelfand--Retakh's quasideterminants, which are defined as follows.\smallskip
 
{\sc Definition} 2.2 [GR]. Let $A=\{\|a_{ij}\|,i\in I,j\in J\}$ be a square 
matrix of order $n={\rm card}\, I={\rm card}\, J$ with formal non-commutative 
entries $a_{ij}$. For $p\in I,q\in J$ denote by $A^{pq}$ the submatrix 
$\{\|a_{ij}\|,i\in I-p,j\in J-q\}$ of $A$. The formula$$|A|_{pq}=a_{pq}-
\sum_{i\in I-p\atop j\in J-q}a_{pj}|A^{pq}|^{-1}_{ij}a_{ip}$$ (which reduces 
to $|A|_{pq}=a_{pq}$ if $n=1$) defines inductively $n^2$ {\it 
quasideterminants} $|A|_{pq}$ of the matrix $A$. (In the commutative case 
$|A|_{pq}=\pm{\rm det}\, A/{\rm det}\, A_{pq}$.)\smallskip 
 
Quasideterminants possess some basic properties of determinants; in 
particular, the following {\it Kramer rule} holds.\smallskip
 
{\sc Proposition} 2.3 [GR]. {\it If $(x_1,\dots ,x_n)$ is the solution of a 
system of equations $$\sum_{j=1}^na_{ij}x_j=\xi_i\, (i=1,\dots ,n),$$ then for 
any $i$$$x_j=|A|_{ij}^{-1}|A_j(\xi)|_{ij},$$ where $A=\|a_{ij}\|$ and 
$A_j(\xi)$ is obtained from $A$ be replacing its $j$-th column by the column 
$(\xi_1,\dots ,\xi_n)$.\smallskip 
 
{\sc Corollary 2.4.} For generic independent $X_1,\dots ,X_n$ and for arbitrary 
$i$\medskip
 
\noindent\quad$-A_j=$ $$\left|\matrix{I&X_1&\dots&X_1^{n-1}\cr 
I&X_2&\dots&X_2^{n-1}\cr \multispan4\quad$\cdots\cdots\cdots\cdots\cdots$\cr 
I&X_n&\dots&X_n^{n-1}\cr}\right|^{-1}_{i,n-j}\left|\matrix{I&X_1&\dots 
&X_1^{n-j-1}&X_1^n&X_1^{n-j+1}&\dots &X_1^{n-1}\cr I&X_2&\dots&X_2^{n-j-
1}&X_2^n&X_2^{n-j+1}&\dots&X_2^{n-1}\cr 
\multispan8\quad$\cdots\cdots\cdots\cdots \cdots\cdots\cdots\cdots 
\cdots\cdots\cdots\cdots \cdots\cdots\cdots$\cr I&X_n&\dots&X_n^{n-j-
1}&X_n^n&X_n^{n-j+1}&\dots&X_n^{n-1}\cr}\right|_{i,n-j}\ (5)$$} 
 
This expression times $(-1)^j$ (with a minor change of notations and with 
$i=n$) is called in [GKLLRT], Section 7.1, the $j$-th elementary symmetric 
function in $X_1,\dots ,X_n$. It is proved in [GKLLRT] and is obvious from 
Corollary 2.4, that for generic $X_1,\dots ,X_n$ it is really symmetric in 
$X_1,\dots ,X_n$.\bigskip 
 
\centerline{\bf 3. The case {\ib n} = 2}\medskip
 
Since for any square matrices $X,Y$$$\left|\matrix{I&I\cr X&Y\cr}\right|={\rm 
det}\, (Y-X),$$the following is precisely Theorem 1.1 for the case 
$n=2$.\smallskip 
 
{\sc Proposition} 3.1. {\it Let $A,B,X,Y$ be square matrices of the same order 
with $X-Y$ being non-degenerate. If$$\eqalign{X^2+AX+B&=0,\cr 
Y^2+AY+B&=0.\cr}\eqno(6)$$Then$$\eqalign{{\rm tr}\, A&=-({\rm tr}\, X+{\rm 
tr}\, Y)\cr {\rm det}\, B&={\rm det}\, X{\rm det}\, Y.\cr}$$} 
 
{\sc Proof.} 1. The equalities (6) imply$$-A=(X^2-Y^2)(X-Y)^{-1},\eqno(7)$$and 
since $X^2-Y^2=(X+Y)(X-Y)+XY-YX$, then$$-A=X+Y+(XY-YX)(X-Y)^{-
1}.$$Therefore$$-{\rm tr}\, A={\rm tr}\, X+{\rm tr}\, Y+{\rm tr}\, XY(X-Y)^{-
1}-{\rm tr}\, YX(X-Y)^{-1},\eqno(8)$$and the relations$$\eqalign{X(X-Y)^{-
1}&=Y(X-Y)^{-1}+I,\cr (X-Y)^{-1}X&=(X-Y)^{-1}Y+I\cr}$$together with the 
identity ${\rm tr}\, UV={\rm tr}\, VU$ imply$$\eqalign{{\rm tr}\, XY(X-Y)^{-
1}&={\rm tr}\, Y(X-Y)^{-1}X={\rm tr}\, Y(X-Y)^{-1}Y+{\rm tr}\, Y,\cr {\rm 
tr}\, YX(X-Y)^{-1}&={\rm tr}\, X(X-Y)^{-1}Y={\rm tr}\, Y(X-Y)^{-1}Y+{\rm tr}\, 
Y.\cr}$$Hence ${\rm tr}\, XY(X-Y)^{-1}={\rm tr}\, YX(X-Y)^{-1}$, and (8) 
yields$$-{\rm tr}\, A={\rm tr}\, X+{\rm tr}\, Y.$$ \smallskip 
 
2. The first of the equalities (6) implies$$B=-(X+A)X,\eqno(9)$$and using (7) 
we get$$\eqalign{(X+A)(X-Y)&=(X-(X^2-Y^2)(X-Y)^{-1})(X-Y)\cr &=X^2-XY-
X^2+Y^2\cr &=Y^2-XY=(X-Y)(-Y).\cr}$$Therefore$$X+A=(X-Y)(-Y)(X-Y)^{-
1},$$$${\rm det}\, (X+A)={\rm det}\, (-Y),$$and in virtue of (9)$${\rm det}\, 
B={\rm det}\, (-(X+A)){\rm det}\, X={\rm det}\, Y{\rm det}\, X.$$\bigskip
 
\centerline{\bf 4. The general case}\medskip
 
Unlike the above proof for $n=2$, our proof in the general case is not reduced 
to a direct calculation and uses the specifics of the matrix 
algebra.\smallskip 
 
{\sc Lemma} 4.1. {\it Let $A_1,\dots ,A_n,X_1,\dots ,X_n$ be square matrices 
of some order $k$, and let the eigenvalues of the matrices $X_1,\dots ,X_n$ be 
$kn$ pairwise different complex numbers. If$$X_i^n+A_1X_i^{n-1}+\dots 
+A_n=0\eqno(10)$$ for $i=1,\dots ,n$, then}$$\eqalign{{\rm tr}\, A_1&=-({\rm 
tr}\, X_1+\dots +{\rm tr}\, X_n),\cr {\rm det}\, A_n&=(-1)^{kn}{\rm det}\, 
X_1\ldots{\rm det}\,X_n.\cr}$$\smallskip 
 
{\sc Proof.} Let $\lambda_{ij},\, j=1,\dots ,k$ be eigenvalues of $X_i$ and 
let $X_iv_{ij}=\lambda_{ij}v_{ij},\, v_{ij}\ne0$. Then for arbitrary $i,j$ 
(10) implies$$\eqalign{0&=(X_i^n+A_1X_i^{n-1}+\dots +A_n)v_{ij}\cr 
&=(\lambda_{ij}^n+\lambda_{ij}^{n-1}A_1+\dots +A_n)v_{ij},\cr}$$whence$${\rm 
det}\,(\lambda_{ij}^n+\lambda_{ij}^{n-1}A_1+\dots +A_n)=0.$$Obviously, 
$$P(\lambda)={\rm det}\,(\lambda^n+\lambda^{n-1}A_1+\dots +A_n)$$is a monic 
polynomial in $\lambda$ of degree $kn$, and since $P(\lambda_{ij})=0$ and all 
$\lambda_{ij}$ are different, then$${\rm det}\, (\lambda^n+\lambda^{n-
1}A_1+\dots +A_n)=\prod_{i,j}(\lambda -\lambda_{ij}).\eqno(11)$$Equate 
constant terms and the coefficients in the term with $\lambda^{kn-1}$ of the 
two sides of (11); we have$$\eqalign{{\rm det}\, A_n&=\prod_{i,j}(-
\lambda_{ij})=(-1)^{kn}\prod_{i,j}\lambda_{ij}=(-1)^{kn}\prod_i{\rm det}\, 
X_i.\cr {\rm tr}A_1&=\sum_{i,j}(-\lambda_{ij})=-\sum_i{\rm tr}\, 
A_i.\cr}$$\smallskip 
 
{\sc Remark 4.2.} Equating other coefficients of polynomials in (11), we may 
get $kn-2$ more identities. Most of them involve matrices, obtained by 
combination of columns of different matrices $A_i$, but in the two extreme 
cases (those of $\lambda^p$ with $p=kn-2$ and 1) we get the formulas which are 
worth mentioning:$${\rm tr}\, A_2+\sigma_2(A_1)=\sum_{1\le i<j\le n}({\rm 
tr}\, X_i{\rm tr}\, X_j)+\sum_{1\le i\le n}\sigma_2(X_i);\eqno(12)$$$${\rm 
det}\, A_{n-1}\cdot{\rm tr}\, A_nA_{n-1}^{-1}=(-1)^{kn-1}(\prod{\rm det}\, 
X_i)(\sum{\rm tr}\, X_i^{-1}),\eqno(13)$$where $\sigma_2$ in (12) denotes the 
second coefficient (the coefficient in the term with $\lambda^{k-2}$) of the 
characteristic polynomial of a matrix, and (13) is valid only if $A_{n-
1},X_1,\dots ,X_n$ are non-degenerate.\smallskip 
 
{\sc Proof of Theorem 1.1.} Let $U,V\subset {\rm Mat}_k{\bb 
C}\times\dots\times{\rm Mat}_k{\bb C}$ ($n$ factors) are respectively the set 
of all $n$-tuples $(X_1,\ldots ,X_n)$ of independent matrices and the set of 
all $n$-tuples $(X_1,\ldots ,X_n)$ of matrices whose eigenvalues are $kn$ 
pairwise different complex numbers. Obviously the both sets are Zariski open 
and non-empty, hence $U\cap V$ is dense in $U$. According to Lemma, the 
equalities (2) both hold in $U\cap V$. Since the both sides of each of the 
equalities (2) are continuous on $U$ (with respect to $X_1,\dots ,X_n$), the 
equalities (2) hold on the whole $U$.\vfill\eject
 
\centerline{\bf 5. A further generalization}\medskip
 
Let $R$ be an associative ring with unity, and $k$ be a field. Suppose that 
there fixed either an additive homomorphism ${\rm tr}\colon R\to k$ satisfying 
the condition ${\rm tr}UV={\rm tr}VU$ for any $U,V\in R$, or a ring 
homomorphism ${\rm det}\colon R\to k$ (or both).\smallskip
 
{\sc Proposition 5.1}. {\it Let $A,B,X,Y\in R$ with $X-Y$ being invertible. If 
$X^2+AX+B=0, Y^2+AY+B=0$, then ${\rm tr}\, A=-({\rm tr}\, X+{\rm tr}\, Y)$ 
and/or ${\rm det}\, B={\rm det}\, X{\rm det}\, Y$, whichever is 
defined.}\smallskip 
 
The proof is the same as in Section 3.
 
To generalize to arbitrary rings the general case of Theorem 1.1 we need an 
explicit expression of $A_1, A_n$ via $X_1,\ldots ,X_n$.\smallskip 
 
{\sc Theorem 5.2}. {\it Let$$\eqalign{A_1&=a_1(X_1,\ldots ,X_n),\cr 
A_n&=a_n(X_1,\ldots ,X_n)\cr}$$be expressions of $A_1,A_n$ via $X_1,\dots ,X_n$ 
involving ring operations and taking inverses and valid where these inverses 
exist} (like (3), (4), (5)). {\it Then for any $X_1,\ldots ,X_n\in R$ each of 
the equalities $$\eqalign{{\rm tr}\, a_1(X_1,\ldots ,X_n)&=-({\rm tr}\, 
X_1+\dots +{\rm tr}\, X_n),\cr {\rm det}\, a_n(X_1\ldots ,X_n)&={\rm det}\, (-
X_1)\ldots{\rm det}\, (-X_n)\cr}$$ holds provided that the both sides exist} 
(that is tr or det is defined and the inverses involved in $a_1$ or $a_n$ 
exist in $R$).\smallskip 
 
Theorem 5.2 cannot be proved by arguments similar to that of Section 4. But it 
is known that identities which hold for complex matrices hold also in 
arbitrary associative rings with unities (see [A], Section 12.4.3). Hence 
Theorem 5.2 follows from Theorem 1.1. \bigskip
Acknowledgement. We are indebted to C. Itzykson and I. Kaplansky for
interesting discussion.  

\centerline{\sc Bibliography}\medskip
 
\item{[A]} {\sc M.~Artin.} Algebra. Prentice Hall, Englewood Cliffs NJ, 
1991.\smallskip
 
\item{[AW]} {\sc S.L.~Adler, Yong-Shi Wu.} Algebraic and geometric aspects of 
generalized quantum dynamics. Preprint hep-th 9405054.\smallskip  
 
\item{[GR]} {\sc I.M.~Gelfand, V.S.~Retakh.} A theory of noncommutative 
determinants and characteristic functions of graphs, I. Publ. du LACIM, Univ. 
de Qu\'ebec \`a Montr\'eal, no. 14, P. 1--26.\smallskip
 
\item{[GKLLRT]} {\sc I.M.~Gelfand, D.~Krob, A.~Lascoux, B.~Leclerc, 
V.S.~Retakh, J.--Y.~Thibon.} Noncommutative symmetric functions. Preprint 
hep-th 9407124.\smallskip 
 
\item{[GS]} {\sc I.M.~Gelfand, M.M.~Smirnov.} The algebra of Chern-Simons 
classes and the Poisson brackets on it. Preprint hep-th 9404103.\smallskip 
 
\item{[K]} {\sc M.~Kontsevich.} Formal (non)-commutative symplectic geometry. 
In: Gelfand Mathematical Seminar, 1992. Birkh\"auser, Boston, 1993, P. 
173--189. 
 
\bye